\documentclass{amsart}
\usepackage{amsfonts,amsmath,amsthm,amssymb}
\newtheorem{theorem}{Theorem}
\newtheorem{corollary}{Corollary}
\newtheorem{proposition}{Proposition}
\title[Integral group ring of Rudvalis simple group]
{Integral group ring of Rudvalis simple group}
\date{}

\author{V.A.~Bovdi, A.B.~Konovalov}

\address{V.A.~Bovdi
\newline Institute of Mathematics, University of Debrecen,
\newline P.O.  Box 12, H-4010 Debrecen, Hungary
\newline Institute of Mathematics and Informatics, College of Ny\'\i regyh\'aza,
\newline S\'ost\'oi \'ut 31/b, H-4410 Ny\'\i regyh\'aza, Hungary}
\email{vbovdi@math.klte.hu}

\address{A.B.~Konovalov
\newline Centre for Interdisciplinary Research in Computational Algebra
\newline School of Computer Science, University of St Andrews,
\newline Jack Cole Building, North Haugh, St Andrews, Fife, KY16 9SX, Scotland}
\email{alexk@mcs.st-andrews.ac.uk}

\subjclass{Primary 16S34, 20C05, secondary 20D08}

\keywords{Zassenhaus conjecture, Kimmerle conjecture,
torsion unit, partial augmentation, integral group ring}

\thanks{Supported by OTKA  No.K68383 and by FAPESP Brasil (proc.08/54650-8)}

\begin{document}
\begin{abstract}
Using the Luthar--Passi method, we investigate the classical
Zassenhaus conjecture for the normalized unit group of the
integral group ring of the Rudvalis sporadic simple group
$\verb"Ru"$. As a consequence, for this group we confirm
Kimmerle's conjecture on prime graphs.
\end{abstract}

\maketitle

\section{Introduction, conjectures and main results}
\label{Intro}

Let $V(\mathbb Z G)$ be  the normalized unit group of the
integral group ring $\mathbb Z G$ of  a finite group $G$. A long-standing
conjecture  of H.~Zassenhaus {\bf (ZC)} says that every torsion unit
$u\in V(\mathbb ZG)$ is conjugate within the rational group algebra
$\mathbb Q G$ to an element in $G$ (see \cite{Zassenhaus}).

For finite simple groups the main tool for the investigation of the
Zassenhaus conjecture is the Luthar--Passi method, introduced in
\cite{Luthar-Passi} to solve it for $A_{5}$ and then applied in 
\cite{Luthar-Trama} for the case of $S_{5}$.
Later M.~Hertweck in \cite{Hertweck1} extended
the Luthar--Passi method and applied it for the investigation
of the Zassenhaus conjecture for $PSL(2,p^{n})$.
The Luthar--Passi method proved to be useful for groups containing
non-trivial normal subgroups as well. For some recent results we refer to
\cite{Bovdi-Hofert-Kimmerle,Bovdi-Konovalov-M11, Hertweck2, Hertweck3,
Hertweck1, Hofert-Kimmerle}.
Also, some related properties and  some weakened variations of the
Zassenhaus  conjecture can be found in
\cite{Artamonov-Bovdi,Bleher-Kimmerle,Luthar-Trama}.

First of all, we need to introduce some notation. By $\# (G)$ we
denote the set of all primes dividing the order of $G$. The
Gruenberg--Kegel graph (or the prime graph) of $G$ is the graph
$\pi (G)$ with vertices labeled by the primes in $\# (G)$ and with
an edge from $p$ to $q$ if there is an element of order $pq$ in
the group $G$. In \cite{Kimmerle} W.~Kimmerle   proposed the
following weakened variation of the Zassenhaus conjecture:

\begin{itemize}
\item[]{\bf (KC)} \qquad
If $G$ is a finite group then $\pi (G) =\pi (V(\mathbb Z G))$.
\end{itemize}

In particular, in the same  paper  W.~Kimmerle verified   that
{\bf (KC)} holds for finite Frobenius and solvable groups. We remark
that with respect to the so-called $p$-version of the Zassenhaus
conjecture the investigation
of Frobenius groups was completed by  M.~Hertweck and the first
author  in \cite{Bovdi-Hertweck}. In \cite{Bovdi-Jespers-Konovalov,
Bovdi-Konovalov-M11, Bovdi-Konovalov-M23, Bovdi-Konovalov-M24, Bovdi-Konovalov-Linton-M22, 
Bovdi-Konovalov-Siciliano-M12} 
{\bf (KC)} was confirmed for the Mathieu simple groups $M_{11}$, 
$M_{12}$, $M_{22}$, $M_{23}$, $M_{24}$ and the sporadic Janko simple groups 
$J_1$, $J_{2}$ and $J_3$.

Here  we continue these investigations for the Rudvalis simple group
$\verb"Ru"$. Although using the Luthar--Passi method we cannot prove the
rational conjugacy for torsion units of $V(\mathbb Z \verb"Ru")$, our
main result gives a lot of information on partial augmentations of
these units. In particular, we confirm the Kimmerle's conjecture for
this group.

Let $G=\verb"Ru"$. It is well known
(see \cite{Ru1}) that
$|G|=2^{14} \cdot 3^{3} \cdot 5^{3} \cdot 7 \cdot 13 \cdot 29$ and
$exp(G)=2^{4} \cdot 3 \cdot 5 \cdot 7 \cdot 13 \cdot 29$.
Let
\[
\begin{split}
\mathcal{C} =\{
C_{1}, C_{2a}, C_{2b}, &  C_{3a},  C_{4a},  C_{4b},  C_{4c},
C_{4d},  C_{5a},  C_{5b},  C_{6a},  C_{7a},  C_{8a},  C_{8b}, \\ C_{8c}, &
C_{10a}, C_{10b},  C_{12a},  C_{12b},  C_{13a},  C_{14a},  C_{14b},  C_{14c},
C_{15a},  C_{16a}, \\ & C_{16b}, C_{20a},  C_{20b},  C_{20c},  C_{24a},  C_{24b},
C_{26a},  C_{26b},  C_{26c},  C_{29a},  C_{29b}
\}
\end{split}
\]
be the collection of all conjugacy classes of $\verb"Ru"$, where the first
index denotes the order of the elements of this conjugacy class
and $C_{1}=\{ 1\}$. Suppose $u=\sum \alpha_g g \in V(\mathbb Z G)$
has finite order $k$. Denote by
$\nu_{nt}=\nu_{nt}(u)=\varepsilon_{C_{nt}}(u)=\sum_{g\in C_{nt}}
\alpha_{g}$ the partial augmentation of $u$ with respect to
$C_{nt}$. 
From the Berman--Higman Theorem 
(see \cite{Berman} and \cite{Sandling}, Ch.5, p.102)
one knows that
$\nu_1 =\alpha_{1}=0$ and
\begin{equation}\label{E:1}
\sum_{C_{nt}\in \mathcal{C}} \nu_{nt}=1.
\end{equation}
Hence, for any character $\chi$ of $G$, we get that $\chi(u)=\sum
\nu_{nt}\chi(h_{nt})$, where $h_{nt}$ is a representative of the
conjugacy class $ C_{nt}$.

Our main result is the following

\begin{theorem}\label{T:1}
Let $G$ denote the Rudvalis sporadic simple group $\verb"Ru"$.
Let $u$ be a torsion unit of $V(\mathbb ZG)$ of order $|u|$ and
let
\[
\begin{split}
\mathfrak{P}(u)=( \nu_{2a},\;&   \nu_{2b},\; \nu_{3a},\;
\nu_{4a},\; \nu_{4b},\; \nu_{4c},\; \nu_{4d},\; \nu_{5a},\;
\nu_{5b},\; \nu_{6a},\; \nu_{7a},\; \nu_{8a},\; \nu_{8b},\\
 \nu_{8c},&\; \nu_{10a},\; \nu_{10b},\; \nu_{12a}, \nu_{12b},\;
\nu_{13a},\; \nu _{14a},\; \nu_{14b},\; \nu_{14c},\; \nu_{15a},\;
\nu_{16a},\\
& \nu_{16b}, \nu_{20a},\; \nu_{20b},\; \nu_{20c},\; \nu_{24a},\;
\nu_{24b},\; \nu_{26a},\; \nu_{26b},\; \nu_{26c},\; \nu_{29a},\;
\nu_{29b}\;) \in \mathbb Z^{35}
\end{split}
\]
be the tuple of partial augmentations of $u$.
The following properties hold.

\begin{itemize}

\item[(i)] If $|u| \not \in 
\{ 28, 30, 40, 48, 52, 56, 60, 80, 104, 112, 120, 208, 240 \}$,
then $|u|$ coincides  with the order of some element $g \in G$.
Equivalently, there is no elements of orders $21$,  $35$,  $39$,
$58$, $65$, $87$, $91$,  $145$,   $203$ and   $377$  in $V(\mathbb
ZG)$.

\item[(ii)] If $|u| \in \{3,7,13\}$, then $u$ is rationally
conjugate to some $g\in G$.

\item[(iii)]
If $|u|=2$, the tuple of the  partial augmentations of $u$ belongs
to the set
\[
\begin{split}
\big\{\;  \mathfrak{P}(u) \; \mid \;
\nu_{2a}+\nu_{2b}=1, \; 
-10 \le \nu_{2a} \le 11, \;
\nu_{kx}=0, \; 
kx \not \in\{2a,2b\}& \; \big\} .
\end{split}
\]

\item[(iv)] If $|u|=5$, the tuple of the  partial augmentations of
$u$ belongs to the set
\[
\begin{split}
\big\{\;  \mathfrak{P}(u) \;  \mid \; 
\nu_{5a}+\nu_{5b}=1, \; 
-1 \le \nu_{5a} \le 6, \;
\nu_{kx}=0,\; 
kx\not\in\{5a,5b\}&\;
\big\} .
\end{split}
\]

\item[(v)] If $|u|=29$, the tuple of the  partial augmentations of
$u$ belongs to the set
\[
\begin{split}
\big\{\; \mathfrak{P}(u) \;  \mid \; 
\nu_{29a}+\nu_{29b}=1, \; 
-4 \le \nu_{29a} \le 5, \;
\nu_{kx}=0,\;
kx\not\in\{29a,29b\}&\;  \big\}.
\end{split}
\]

\end{itemize}
\end{theorem}

As an immediate consequence of  part (i) of the Theorem we obtain

\begin{corollary} If $G = \verb"Ru"$ then
$\pi(G)=\pi(V(\mathbb ZG))$.
\end{corollary}

\section{Preliminaries}
The following result is a reformulation of
the Zassenhaus conjecture in terms of vanishing 
of partial augmentations of torsion units.

\begin{proposition}\label{P:5}
(see \cite{Luthar-Passi} and
Theorem 2.5 in \cite{Marciniak-Ritter-Sehgal-Weiss})
Let $u\in V(\mathbb Z G)$
be of order $k$. Then $u$ is conjugate in $\mathbb
QG$ to an element $g \in G$ if and only if for
each $d$ dividing $k$ there is precisely one
conjugacy class $C$ with partial augmentation
$\varepsilon_{C}(u^d) \neq 0 $.
\end{proposition}

The next result now yield that several partial augmentations
are zero.

\begin{proposition}\label{P:4}
(see \cite{Hertweck2}, Proposition 3.1;
\cite{Hertweck1}, Proposition 2.2)
Let $G$ be a finite
group and let $u$ be a torsion unit in $V(\mathbb
ZG)$. If $x$ is an element of $G$ whose $p$-part,
for some prime $p$, has order strictly greater
than the order of the $p$-part of $u$, then
$\varepsilon_x(u)=0$.
\end{proposition}

The key restriction on partial augmentations is given
by the following result that is the cornerstone of
the Luthar--Passi method.

\begin{proposition}\label{P:1}
(see \cite{Luthar-Passi,Hertweck1}) Let either $p=0$ or $p$ a prime
divisor of $|G|$. Suppose
that $u\in V( \mathbb Z G) $ has finite order $k$ and assume $k$ and
$p$ are coprime in case $p\neq 0$. If $z$ is a complex primitive $k$-th root
of unity and $\chi$ is either a classical character or a $p$-Brauer
character of $G$, then for every integer $l$ the number
\begin{equation}\label{E:2}
\mu_l(u,\chi, p ) =
\textstyle\frac{1}{k} \sum_{d|k}Tr_{ \mathbb Q (z^d)/ \mathbb Q }
\{\chi(u^d)z^{-dl}\}
\end{equation}
is a non-negative integer.
\end{proposition}

Note that if $p=0$, we will use the notation $\mu_l(u,\chi
, * )$ for $\mu_l(u,\chi , 0)$.

Finally, we shall use the well-known bound for
orders of torsion units.

\begin{proposition}\label{P:2}  (see  \cite{Cohn-Livingstone})
The order of a torsion element $u\in V(\mathbb ZG)$
is a divisor of the exponent of $G$.
\end{proposition}

\section{Proof of the Theorem}

Throughout this section we denote the group $\verb"Ru"$ by $G$.
The character table of $G$,
as well as the $p$-Brauer character tables, which will be denoted by
$\mathfrak{BCT}{(p)}$ where $p\in\{2,3,5,7,13,29\}$, can be found using
the computational algebra system GAP \cite{GAP},
which derives these data from \cite{AFG,ABC}. For the characters
and conjugacy classes we will use throughout the paper the same
notation, indexation inclusive, as used in the GAP Character Table Library.

Since the group $G$ possesses elements of orders $2$, $3$, $4$,
$5$, $6$, $7$, $8$, $10$,  $12$,  $13$,  $14$,  $15$,  $16$, $20$,
$24$, $26$ and  $29$, first of all we  investigate  units of some
these orders (except the units of orders $4$,  $6$,  $8$,  $10$,
$12$, $14$, $15$, $16$, $20$,  $24$ and   $26$). After this, by
Proposition \ref{P:2}, the order of each torsion unit divides the
exponent of $G$, so to prove the Kimmerle's conjecture, it remains to 
consider units of orders $21$, $35$, $39$,
$58$, $65$, $87$, $91$, $145$, $203$ and $377$. We will prove 
that no units of all these orders do appear in $V(\mathbb ZG)$.

Now we consider each case separately.

\noindent$\bullet$ Let $u$ be an involution. By (\ref{E:1}) and
Proposition \ref{P:4} we have that $\nu_{2a}+\nu_{2b}=1$. Put $t_1
= 3 \nu_{2a} - 7 \nu_{2b}$ and $t_2=11\nu_{2a} - 7 \nu_{2b}$.
Applying Proposition \ref{P:1} we get the following system
\[
\begin{split}
\mu_{1}(u,\chi_{2},*) & = \textstyle \frac{1}{2} (2t_1 + 378) \geq
0; \qquad
\mu_{0}(u,\chi_{2},*)  = \textstyle \frac{1}{2} (-2t_1 + 378) \geq 0; \\ 
\mu_{0}(u,\chi_{4},*) & = \textstyle \frac{1}{2} (2t_2 + 406) \geq
0; \qquad
\mu_{1}(u,\chi_{4},*) = \textstyle \frac{1}{2} (-2t_2 + 406) \geq 0. \\ 
\end{split}
\]
From these restrictions and the requirement that all
$\mu_i(u,\chi_{j},*)$ must be non-ne\-ga\-tive integers we get
$22$ pairs $(\nu_{2a},\nu_{2b})$ listed in part (iii) of the
Theorem \ref{T:1}.

Note that using our implementation of the Luthar--Passi method, 
which we intended to make available in the GAP package LAGUNA \cite{LAGUNA},
we computed inequalities from Proposition \ref{P:1} for every
irreducible character from ordinary and Brauer character tables,
and for every $0 \leq l \leq |u|-1$, but the only inequalities that really
matter are those ones listed above. The same remark applies for
all other orders of torsion units considered in the paper.

\noindent $\bullet$ Let $u$ be a unit of order either $3$, $7$ or
$13$. Using Proposition \ref{P:4} we obtain that all partial
augmentations except one are zero. Thus by Proposition \ref{P:5}
part (ii) of the Theorem \ref{T:1} is proved.

\noindent $\bullet$ Let $u$ be a unit of order $5$. By (\ref{E:1})
and Proposition \ref{P:4} we get $\nu_{5a}+\nu_{5b}=1$. Put
$t_1 = 6 \nu_{5a} +  \nu_{5b}$ and $t_2=3 \nu_{5a} - 2 \nu_{5b}$.
By (\ref{E:2}) we obtain the system of inequalities
\[
\begin{split}
\mu_{0}(u,\chi_{4},*) & = \textstyle \frac{1}{5} (4t_1 + 406) \geq
0; \qquad \mu_{1}(u,\chi_{4},*)  = \textstyle \frac{1}{5} (-t_1+ 406) \geq 0; \\ 
\mu_{0}(u,\chi_{2},2) & = \textstyle \frac{1}{5} (4t_2 + 28) \geq
0; \qquad  \mu_{1}(u,\chi_{2},2)  = \textstyle \frac{1}{5} (-t_2 + 28) \geq 0. \\ 
\end{split}
\]
Again, using the condition for $\mu_i(u,\chi_{j},p)$
to be non-negative integers, we obtain eight pairs
$(\nu_{5a},\nu_{5b})$ listed in part (iv) of the Theorem
\ref{T:1}.

\noindent$\bullet$ Let $u$ be a unit of order $29$. By (\ref{E:1})
and Proposition \ref{P:4} we have that $\nu_{29a}+\nu_{29b}=1$.
Put $t_1 = 15 \nu_{29a} - 14 \nu_{29b}$. Then using (\ref{E:2})
we obtain the system of inequalities
\[
\begin{split}
\mu_{1}(u,\chi_{6},2) & = \textstyle \frac{1}{29} (t_1 + 8192)
\geq 0; \quad \mu_{2}(u,\chi_{7},5)  = \textstyle \frac{1}{29} (-t_1 + 2219) \geq 0; \\ 
&\mu_{1}(u,\chi_{2},5)  = \textstyle \frac{1}{29} (12 \nu_{29a} - 17 \nu_{29b} + 133) \geq 0; \\ 
&\mu_{2}(u,\chi_{2},5)  = \textstyle \frac{1}{29} (-17 \nu_{29a} + 12 \nu_{29b} + 133) \geq 0. \\ 
\end{split}
\]
Now applying the condition for $\mu_i(u,\chi_{j},p)$ to be
non-negative integers we obtain ten  pairs
$(\nu_{29a},\nu_{29b})$ listed in part (v) of the Theorem
\ref{T:1}.

Now it remains to prove part (i) of the Theorem \ref{T:1}.

\noindent$\bullet$ Let $u$ be a unit of order $21$. By (\ref{E:1})
and Proposition \ref{P:4} we obtain that $\nu_{3a}+\nu_{7a}=1$.
By (\ref{E:2}) we obtain the system of inequalities
\[
\begin{split}
\mu_{1}(u,\chi_{4},*) & = \textstyle \frac{1}{21} ( \nu_{3a} +
405) \geq 0; \qquad
\mu_{0}(u,\chi_{2},2)  = \textstyle \frac{1}{21} (12 \nu_{3a} + 30) \geq 0; \\ 
&\qquad \mu_{7}(u,\chi_{2},2)  = \textstyle \frac{1}{21} (-6 \nu_{3a} + 27) \geq 0, \\ 
\end{split}
\]
which has no integer solutions such that all $\mu_i(u,\chi_{j},p)$
are non-negative integers.

\noindent$\bullet$ Let $u$ be a unit of order $35$. By (\ref{E:1})
and Proposition \ref{P:4} we get $\nu_{5a}+\nu_{7a}+\nu_{7b}=1$.
Put $t_1 =  \nu_{5a} +  \nu_{5b}$. Since $|u^{7}|=5$, for any
character $\chi$ of $G$ we need to consider eight cases defined by
part (iv) of the Theorem. Using (\ref{E:2}), in all of these
cases we get the same system of inequalities
\[
\begin{split}
\mu_{0}(u,\chi_{2},*) & = \textstyle \frac{1}{35} (72 t_1 + 390) \geq 0; \qquad
\mu_{0}(u,\chi_{4},2)  = \textstyle \frac{1}{35} (-96 t_1+ 1230) \geq 0, \\ 
\end{split}
\]
which has no integer solutions such that all $\mu_i(u,\chi_{j},p)$
are non-negative integers.

\noindent$\bullet$ Let $u$ be a unit of order $39$. By (\ref{E:1})
and Proposition \ref{P:4} we have that $\nu_{3a}+\nu_{13a}=1$.
By (\ref{E:2}) we obtain that
\[
\begin{split}
\mu_{0}(u,\chi_{5},*) & = \textstyle \frac{1}{39} (72 \nu_{13a} + 819) \geq 0; \quad
\mu_{13}(u,\chi_{5},*)  = \textstyle \frac{1}{39} (- 36 \nu_{13a} + 819) \geq 0; \\ 
\mu_{1}(u,\chi_{2},*) & = \textstyle \frac{1}{39} (\nu_{13a} + 377) \geq 0; \quad
\mu_{13}(u,\chi_{2},2)  = \textstyle \frac{1}{39} (-12 \nu_{3a} - 24 \nu_{13a} + 51) \geq 0. \\ 
\end{split}
\]
From the first two inequalities we obtain that $\nu_{13a} \in \{0,13\}$,
and now the last two inequalities lead us to a contradiction.

\noindent$\bullet$ Let $u$ be a unit of order $58$. By (\ref{E:1})
and Proposition \ref{P:4} we have that  \
$$
\nu_{2a}+\nu_{2b}+\nu_{29a}+\nu_{29b}=1.
$$
Put $t_1 = 6 \nu_{2a} - 14 \nu_{2b} -  \nu_{29a} -  \nu_{29b}$,
$t_2 = 11\nu_{2a} - 7\nu_{2b}$ and $t_3 =64 \nu_{2b} + 14
\nu_{29a} - 15 \nu_{29b}$. Since $|u^2|=29$ and $|u^{29}|=2$,
according to parts (iii) and (v) of the Theorem we need to
consider 220 cases, which we can group in the following way.
First, let
\[
\begin{split}
\chi(u^{29}) \in
\{ \; &
\chi(2a), \;
-5\chi(2a)+6\chi(2b), \;
-10\chi(2a)+11\chi(2b), \\
& -2\chi(2a)+3\chi(2b), \;
-8\chi(2a)+9\chi(2b), \;
6\chi(2a)-5\chi(2b), \\
& 3\chi(2a)-2\chi(2b), \;
9\chi(2a)-8\chi(2b), \;
4\chi(2a)-3\chi(2b) \; \}.
\end{split}
\]
 Then
by (\ref{E:2}) we obtain the system of inequalities
\[
\begin{split}
\mu_{0}(u,\chi_{2},*) = \textstyle \frac{1}{58} & (-28 t_1 + \alpha) \geq 0; \quad 
\mu_{29}(u,\chi_{2},*) = \textstyle \frac{1}{58} (28 t_1 + \beta) \geq 0; \\ 
&\mu_{1}(u,\chi_{2},*) = \textstyle \frac{1}{58} (- t_1 + \gamma) \geq 0, \\ 
\end{split}
\]
$$ \text{where} \quad (\alpha,\beta,\gamma) = {\tiny
\begin{cases}
(400,412,383), \; \text{if} \; \chi(u^{29}) = \chi(2a) ;\\
(520,292,263), \; \text{if} \; \chi(u^{29}) = -5\chi(2a)+6\chi(2b) ;\\
(620,192,163), \; \text{if} \; \chi(u^{29}) = -10\chi(2a)+11\chi(2b); \\
(460,352,323), \; \text{if} \; \chi(u^{29}) = -2\chi(2a)+3\chi(2b); \\
(580,232,203), \; \text{if} \; \chi(u^{29}) = -8\chi(2a)+9\chi(2b); \\
(300,512,483), \; \text{if} \; \chi(u^{29}) = 6\chi(2a)-5\chi(2b) ; \\
(360,452,423), \; \text{if} \; \chi(u^{29}) = 3\chi(2a)-2\chi(2b) ; \\
(240,572,543), \; \text{if} \; \chi(u^{29}) = 9\chi(2a)-8\chi(2b) ; \\
(340,472,443), \; \text{if} \; \chi(u^{29}) = 4\chi(2a)-3\chi(2b), \\
\end{cases}}
$$ which  has no integral solution such that all
$\mu_{i}(u,\chi_{j},p)$ are non-negative integers.

In the remaining cases we consider the following system obtained by (\ref{E:2}):
\[
\begin{split}
\mu_{0}(u,\chi_{2},*) & = \textstyle \frac{1}{58} (-28t_1 +
\alpha_1) \geq 0; \qquad
\mu_{29}(u,\chi_{2},*)  = \textstyle \frac{1}{58} (28t_1 + \alpha_2) \geq 0; \\ 
\mu_{0}(u,\chi_{4},*) & = \textstyle \frac{1}{58} (56t_2 +\alpha_3) \geq 0;\quad  \qquad \mu_{29}(u,\chi_{4},*) = \textstyle \frac{1}{58} (-56t_2+ \alpha_4) \geq 0; \\ 
\mu_{1}(u,\chi_{34},*) & = \textstyle \frac{1}{58} (-t_3+ \beta_1)
\geq 0; \qquad\quad
\mu_{4}(u,\chi_{34},*)  = \textstyle \frac{1}{58} (t_3 + \beta_2) \geq 0, \\ 
\end{split}
\]
where the tuple of coefficients $(\alpha_1,\alpha_2,\alpha_3,\alpha_4)$
depends only of the value of $\chi(u^{29})$:
$$
(\alpha_1,\alpha_2,\alpha_3,\alpha_4) =
{\tiny
\begin{cases}
(420, 392, 392, 420), & \; \text{if} \; \chi(u^{29}) = \chi(2b);\\
(320, 492, 572, 240), & \; \text{if} \; \chi(u^{29}) = 5\chi(2a)-4\chi(2b);\\
(600, 212, 68, 744),  & \; \text{if} \; \chi(u^{29}) = -9\chi(2a)+10\chi(2b);\\
(540, 272, 176, 636), & \; \text{if} \; \chi(u^{29}) = -6\chi(2a)+7\chi(2b);\\
(380, 432, 464, 348), & \; \text{if} \; \chi(u^{29}) = 2\chi(2a)-\chi(2b);\\
(260, 552, 680, 132), & \; \text{if} \; \chi(u^{29}) = 8\chi(2a)-7\chi(2b);\\ 
(480, 332, 284, 528), & \; \text{if} \; \chi(u^{29}) = -3\chi(2a)+4\chi(2b);\\ 
(500, 312, 248, 564), & \; \text{if} \; \chi(u^{29}) = -4\chi(2a)+5\chi(2b);\\ 
(200, 612, 24, 24),  & \; \text{if} \; \chi(u^{29}) = 11\chi(2a)-10\chi(2b);\\ 
(220, 592, 752, 60), & \; \text{if} \; \chi(u^{29}) = 10\chi(2a)-9\chi(2b);\\ 
(280, 532, 644, 168), & \; \text{if} \; \chi(u^{29}) = 7\chi(2a)-6\chi(2b);\\ 
(440, 372, 356, 456), & \; \text{if} \; \chi(u^{29}) = -\chi(2a)+2\chi(2b);\\ 
(560, 252, 140, 672), & \; \text{if} \; \chi(u^{29}) = -7\chi(2a)+8\chi(2b),\\ 
\end{cases}}
$$
while the pair $(\beta_1,\beta_2)$ depends both on $\chi(u^{29})$ and $\chi(u^{2})$:
$$ \tiny{ \begin{array}{|c|c|c|c|c|} \hline
& \chi(29a) 
& \chi(29b) 
& 5\chi(29a)-4\chi(29b) 
& -2\chi(29a)+3\chi(29b) \\ \hline 
\chi(2b) & 110641, 110513 & 110670, 110542 & 110525, 110397 & 110728, 110600 \\ 
5\chi(2a)-4\chi(2b) & 110321, 110833 & 110350, 110862 & 110205, 110717 & 110408, 110920 \\ 
-9\chi(2a)+10\chi(2b) & 111217, 109937 & 111246, 109966 &  111101, 109821 & 111304, 110024 \\ 
-6\chi(2a)+7\chi(2b) & 111025, 110129 & 111054, 110158 & 110909, 110013 & 111112, 110216 \\ 
2\chi(2a)-\chi(2b) & 110513, 110641 &  110542, 110670 & 110397, 110525 & 110600, 110728 \\ 
8\chi(2a)-7\chi(2b) & 110129, 111025 & 110158, 111054 & 110013, 110909 & 110216, 111112 \\ 
-3\chi(2a)+4\chi(2b) & 110833, 110321 & 110862, 110350 & 110717, 110205 & 110920, 110408 \\ 
-4\chi(2a)+5\chi(2b) & 110897, 110257 &  110926, 110286 & 110781, 110141 & 110984, 110344 \\ 
11\chi(2a)-10\chi(2b) & 109937, 111217 & 109966, 111246 & 109821, 111101 & 110024, 111304 \\ 
10\chi(2a)-9\chi(2b) & 110001, 111153 & 110030, 111182 & 109885, 111037 & 110088, 111240 \\ 
7\chi(2a)-6\chi(2b) & 110193, 110961 & 110222, 110990 & 110077, 110845 & 110280, 111048\\ 
-\chi(2a)+2\chi(2b) & 110705, 110449 & 110734, 110478 & 110589, 110333 & 110792, 110536  \\ 
-7\chi(2a)+8\chi(2b) & 111089, 110065 & 111118, 110094 & 110973, 109949 & 111176, 110152 \\ 
\hline \end{array}} $$
$$ \tiny{ \begin{array}{|c|c|c|c|} \hline
& 2\chi(29a)-\chi(29b) 
& -3\chi(29a)+4\chi(29b) 
& -4\chi(29a)+5\chi(29b) \\ \hline 
\chi(2b) & 110612, 110484 & 110757, 110629 & 110786, 110658 \\ 
5\chi(2a)-4\chi(2b) & 110292, 110804 & 110437, 110949 &  110466, 110978 \\ 
-9\chi(2a)+10\chi(2b) & 111188, 109908 & 111333, 110053 &  111362, 110082  \\ 
-6\chi(2a)+7\chi(2b) & 110996, 110100 & 111141, 110245 & 111170, 110274 \\ 
2\chi(2a)-\chi(2b) & 110484, 110612 & 110629, 110757 & 110658, 110786 \\ 
8\chi(2a)-7\chi(2b) & 110100, 110996 & 110245, 111141 & 110274, 111170 \\ 
-3\chi(2a)+4\chi(2b) &  110804, 110292 & 110949, 110437 & 110978, 110466 \\ 
-4\chi(2a)+5\chi(2b) & 110868, 110228 & 111013, 110373 & 111042, 110402 \\ 
11\chi(2a)-10\chi(2b) & 109908, 111188 & 110053, 111333 & 110082, 111362 \\ 
10\chi(2a)-9\chi(2b) & 109972, 111124 & 110117, 111269 & 110146, 111298 \\ 
7\chi(2a)-6\chi(2b) & 110164, 110932 & 110309, 111077 & 110338, 111106 \\ 
-\chi(2a)+2\chi(2b) & 110676, 110420 & 110821, 110565 & 110850, 110594 \\ 
-7\chi(2a)+8\chi(2b) & 111060, 110036 & 111205, 110181 & 111234, 110210 \\ 
\hline 
\hline
& 3\chi(29a)-2\chi(29b) 
& -\chi(29a)+2\chi(29b) 
& 4\chi(29a)-3\chi(29b) \\ \hline 
\chi(2b) & 110583, 110455 & 110699, 110571 & 110554, 110426 \\ 
5\chi(2a)-4\chi(2b) & 110263, 110775 & 110379, 110891 & 110234, 110746 \\ 
-9\chi(2a)+10\chi(2b) &  111159, 109879 & 111275, 109995 & 111130, 109850 \\ 
-6\chi(2a)+7\chi(2b) & 110967, 110071 & 111083, 110187 & 110938, 110042 \\ 
2\chi(2a)-\chi(2b) & 110455, 110583 & 110571, 110699 & 110426, 110554 \\ 
8\chi(2a)-7\chi(2b) & 110071, 110967 & 110187, 11108 & 110042, 110938 \\ 
-3\chi(2a)+4\chi(2b) & 110775, 110263 & 110891, 110379 & 110746, 110234 \\ 
-4\chi(2a)+5\chi(2b) & 110839, 110199 & 110955, 110315 &  110810, 110170 \\ 
11\chi(2a)-10\chi(2b) & 109879, 111159 &  109995, 111275 & 109850, 111130 \\ 
10\chi(2a)-9\chi(2b) & 109943, 111095 & 110059, 111211 & 109914, 111066 \\ 
7\chi(2a)-6\chi(2b) & 110135, 110903 &  110251, 111019 & 110106, 110874 \\ 
-\chi(2a)+2\chi(2b) & 110647, 110391 & 110763, 110507 &  110618, 110362 \\ 
-7\chi(2a)+8\chi(2b) & 111031, 110007 & 111147, 110123 & 111002, 109978 \\ 
\hline \end{array}} $$
Additionally, when $\chi(u^{29}) \in \{ \chi(2b), 7\chi(2a)-6\chi(2b), -7\chi(2a)+8\chi(2b) \}$,
we need to consider one more inequality
$$
\mu_{1}(u,\chi_{2},*) = \textstyle \frac{1}{58} (-6 \nu_{2a} + 14 \nu_{2b} +  \nu_{29a} +  \nu_{29b} + \gamma) \geq 0, \\ 
$$
where\qquad  $ \gamma = {\tiny
\begin{cases}
363, & \; \text{if} \; \chi(u^{29}) = \chi(2b);\\
503, & \; \text{if} \; \chi(u^{29}) = 7\chi(2a)-6\chi(2b);\\ 
223, & \; \text{if} \; \chi(u^{29}) = -7\chi(2a)+8\chi(2b).\\ 
\end{cases}}
$ \noindent \newline All systems of inequalities, constructed as
described above, have no integer solutions such that all
$\mu_{i}(u,\chi_{j},p)$ are non-negative integers.

\noindent$\bullet$ Let $u$ be a unit of order $65$.
By (\ref{E:1}) and Proposition \ref{P:4} we have that
$$
\nu_{5a}+\nu_{5b}+\nu_{13a}=1.
$$
Since $|u^{13}|=5$, we need to consider eight cases
listed in part (iv) of the Theorem.
Put
$t_1 = 3 \nu_{5a} + 3 \nu_{5b} +  \nu_{13a}$
and
$t_2 = 6 \nu_{5a} +  \nu_{5b} + 3 \nu_{13a}$.
Then using (\ref{E:2}) we obtain
\[
\begin{split}
\mu_{0}(u,\chi_{2},*) & = \textstyle \frac{1}{65} (48t_1 + 402) \geq 0; \qquad
\mu_{13}(u,\chi_{2},*)  = \textstyle \frac{1}{65} (-12t_1 + 387) \geq 0; \\ 
\mu_{0}(u,\chi_{4},*) & = \textstyle \frac{1}{65} (48t_2 + \alpha ) \geq 0; \qquad
\mu_{13}(u,\chi_{4},*) = \textstyle \frac{1}{65} (-12 t_2 + \beta ) \geq 0, \\ 
\end{split}
\]
$$ \text{where} \quad (\alpha,\beta) =
{\tiny
\begin{cases}
(466,436), \; \text{if} \; \chi(u^{13}) = \chi(5a) ;\\
(446,441), \; \text{if} \; \chi(u^{13}) = \chi(5b) ;\\
(546,416), \; \text{if} \; \chi(u^{13}) = 5\chi(5a)-4\chi(5b) ;\\
(486,431), \; \text{if} \; \chi(u^{13}) = 2\chi(5a)-\chi(5b) ;\\
(566,411), \; \text{if} \; \chi(u^{13}) = 6\chi(5a)-5\chi(5b) ;\\
(506,426), \; \text{if} \; \chi(u^{13}) = 3\chi(5a)-2\chi(5b) ;\\
(426,446), \; \text{if} \; \chi(u^{13}) = -\chi(5a)+2\chi(5b) ;\\
(526,421), \; \text{if} \; \chi(u^{13}) = 4\chi(5a)-3\chi(5b).
\end{cases}}$$

\noindent In all cases we have no solutions such that all
$\mu_{i}(u,\chi_{i},p)$ are non-negative integers.

\noindent$\bullet$ Let $u$ be a unit of order $87$. By (\ref{E:1})
and Proposition \ref{P:4} we have that
$$
\nu_{3a}+\nu_{29a}+\nu_{29b}=1.
$$
Since $|u^3|=29$, according to part (v) of the Theorem we need to
consider ten cases.  Put $t_1 = \nu_{29a} + \nu_{29b}$. In all of
these cases by (\ref{E:2}) we get the system
\[
\begin{split}
\mu_{0}(u,\chi_{2},*) & = \textstyle \frac{1}{87} (56t_1 + 406)
\geq 0; \qquad
\mu_{29}(u,\chi_{2},*) = \textstyle \frac{1}{87} (- 28t_1 + 406) \geq 0, \\ 
\end{split}
\]
that lead us to a contradiction.

\noindent$\bullet$ Let $u$ be a unit of order $91$. By (\ref{E:1})
and Proposition \ref{P:4} we get $ \nu_{7a}+\nu_{13a}=1$.
Now using (\ref{E:2}) we obtain non-compatible inequalities
\[
\begin{split}
\mu_{0}(u,\chi_{2},2) & = \textstyle \frac{1}{91} (144 \nu_{13a} +
52) \geq 0; \qquad
\mu_{7}(u,\chi_{2},2)  = \textstyle \frac{1}{91} (- 12 \nu_{13a} + 26) \geq 0. \\ 
\end{split}
\]

\noindent$\bullet$ Let $u$ be a unit of order $145$. By
(\ref{E:1}) and Proposition \ref{P:4} we have that
$$
\nu_{5a}+\nu_{5b}+\nu_{29a}+\nu_{29b}=1.
$$
Put $t_1 = 3 \nu_{5a} + 3 \nu_{5b} +  \nu_{29a} +  \nu_{29b}$.
Since $|u^{29}|=5$ and $|u^{5}|=29$, for any character $\chi$ of $G$
we need to consider $80$ cases defined by parts
(iv) and (v) of the Theorem. Luckily, in every case by (\ref{E:2})
we obtain the same pair of incompatible inequalities
\[
\begin{split}
\mu_{0}(u,\chi_{2},*) & = \textstyle \frac{1}{145} (112 t_1 + 418)
\geq 0; \qquad
\mu_{29}(u,\chi_{2},*)  = \textstyle \frac{1}{145} (-28t_1 + 403) \geq 0. \\ 
\end{split}
\]

\noindent$\bullet$ Let $u$ be a unit of order $203$. By
(\ref{E:1}) and Proposition \ref{P:4} we have that
$$
\nu_{7a}+\nu_{29a}+\nu_{29b}=1.
$$
Since $|u^7|=29$, according to part (v) of the Theorem
we need to consider ten cases.
Put $t_1 = \nu_{29a} +  \nu_{29b}$, and then using (\ref{E:2})
in each case we obtain a non-compatible system of inequalities
\[
\begin{split}
\mu_{29}(u,\chi_{2},2) & = \textstyle \frac{1}{203} (28 t_1) \geq 0; 
\qquad \mu_{0}(u,\chi_{2},2)  = \textstyle \frac{1}{203} (- 168 t_1) \geq 0; \\ 
&\mu_{1}(u,\chi_{2},*)  = \textstyle \frac{1}{203} (t_1 + 377) \geq 0. \\ 
\end{split}
\]

\noindent$\bullet$ Let $u$ be a unit of order $377$. By
(\ref{E:1}) and Proposition \ref{P:4} we have that
$$
\nu_{13a}+\nu_{29a}+\nu_{29b}=1.
$$
Since $|u^{13}|=29$,  we need to consider ten cases defined by
part (v) of the Theorem. In each case by (\ref{E:2}) we obtain
the following system of inequalities
\[
\begin{split}
\mu_{0}(u,\chi_{4},*) & = \textstyle \frac{1}{377} (1008 \nu_{13a}
+ 442) \geq 0; \\
\mu_{29}(u,\chi_{4},*) &= \textstyle \frac{1}{377} (-84 \nu_{13a} + 403) \geq 0. \\ 
\end{split}
\]
which have no solution such that all $\mu_{i}(u,\chi_{j},*)$ are non-negative integers.

\bibliographystyle{plain}
\bibliography{Bovdi_Konovalov_Ru}

\end{document}